\documentclass[reqno,12pt]{amsart}
\theoremstyle{plain}    
\newtheorem{thm}{Theorem}[section]
\newtheorem{lemma}[thm]{Lemma} 
\newtheorem{prop}[thm]{Proposition}

\theoremstyle{remark}

\newtheorem{remark}[thm]{Remark}

\theoremstyle{definition}

\newtheorem{notation}[thm]{Notation}

\usepackage{amscd,amssymb,comment,epic,eepic,enumerate,euscript,graphics}
\usepackage[initials]{amsrefs}

\newcommand\Ac{{\mathcal{A}}}
\newcommand\At{{\widetilde A}}
\newcommand\Bc{{\mathcal{B}}}
\newcommand\Cpx{{\mathbf C}}
\newcommand\ct{{\tilde c}}
\newcommand\Fb{{\mathbf F}}
\newcommand\fdim{\text{\rm fdim}\,}
\newcommand\Ints{{\mathbf Z}}
\newcommand\Lambdao{{\Lambda\oup}}
\newcommand\lspan{\mathrm{span}\,}
\newcommand\Mcal{{\mathcal{M}}} 
\newcommand\Mcalt{{\widetilde\Mcal}}
\newcommand\Nats{{\mathbf N}}
\newcommand\Nc{{\mathcal{N}}}
\newcommand\oup{^{\mathrm o}}
\newcommand\restrict{{\upharpoonright}}
\newcommand\REu{{\EuScript R}}                   
\newcommand\SEu{{\EuScript S}}
\newcommand\smd[2]{\underset{#2}{#1}}
\newcommand\smdp[3]{\overset{#3}{\smd{#1}{#2}}}

\begin{document}

\title[amalgamated free products]{A description of amalgamated free products of finite von Neumann algebras
over finite dimensional subalgebras}

\author[Dykema]{Ken Dykema$^{*}$}
\address{Department of Mathematics, Texas A\&M University,
College Station, TX 77843-3368, USA}
\email{kdykema@math.tamu.edu}
\thanks{\footnotesize $^{*}$Research supported in part by NSF grant DMS-0901220}



\date{February 10, 2010}

\begin{abstract}
We show that a free product of a II$_1$--factor and a finite
von Neumann algebra with amalgamation over a finite dimensional subalgebra
is always a II$_1$--factor, and provide
an algorithm for describing it in terms of free products (with amalgamation over the scalars)
and compression/dilation.
As an application, we show that the class of direct sums of finitely many von Neumann algebras that are
interpolated free group factors, hyperfinite II$_1$--factors, type I$_n$ algebras for $n$ finite, and
finite dimensional algebras, is closed under taking free products with amalgamation over finite dimensional
subalgebras.
\end{abstract}

\maketitle

\section{Introduction}

The amalgamated free product of von Neumann algebras has been an important construction in subfactor theory, since the work~\cite{Po93}
of S.\ Popa,
where this construction was used in producing subfactors with arbitrary (allowable) index.
Other (but more limited) constructions of subfactors involving free products appeared in~\cite{Ra94} and~\cite{D95}.
Amalgamated free products appeared again in~\cite{Po95}, where Popa constructed subfactors with arbitrary (allowable) higher relative commutants.
More recently, in~\cite{GJS1}, Guionnet, Jones and Shlyakhtenko gave another construction of such subfactors, involving planar algebras.

In this note, we consider free products 
\begin{equation}\label{eq:ADBintro}
\Mcal=A*_DB
\end{equation}
of certain sorts of finite von Neumann algebras $A$ and $B$ with amalgamation over a finite dimensional subalgebra $D$,
and describe how to obtain $\Mcal$ by taking free products (over $\Cpx$) and operations of rescaling.
These results are used in other authors' investigations of and results on the subfactors constructed in~\cite{GJS1}.

Our main result (Theorem~\ref{thm:M=ADB}) is a precise statement of the algorithm for describing the algebra $\Mcal$ in the free
product~\eqref{eq:ADBintro} if one of $A$ or $B$ is a II$_1$--factor.
Previous results, found in~\cite{D93}, \cite{D95} and~\cite{BD04}, expressed respectively (i) free products of finite dimensional and hyperfinite
von Neumann algebras (with amalgamation over the scalars $\Cpx$),
(ii) amagamated free products of finite dimensional von Neumann algebras and, more generally, type~I von Neumann algebras
having atomic centers
and (iii) the amalgamated free product of two copies of the hyperfinite II$_1$--factors over a type~I subalgebra
with atomic center,
in terms of interpolated free group factors $L(\Fb_t)$, $1<t\le\infty$ (see~\cite{Ra94}, \cite{D94:i})
and hyperfinite von Neumann algebras.

As an application of our main result, we consider 
the case when $A$ and $B$ are both in the class of direct sums of finitely many subalgebras which may be finite type I algebras,
hyperfinite II$_1$--factors and interpolated free group factors.
We show that this class is closed under taking free products with amalgamation over finite dimensional subalgebras,
and, moreover, that if all the interpolated free group factors $L(\Fb_t)$ appearing as direct summands in $A$ and $B$
have finite parameter $t$, then so do all interpolated free group factors appearing as direct summands in $\Mcal$.
This latter result is used by Kodiyalam and Sunder in~\cite{KS}
and also in the related paper~\cite{GJS} of Guionnet, Jones and Shlyakhtenko.

\smallskip
\noindent
{\em Acknowledgement:}
This work was inspired by a question from
Vijay Kodiyalam and V.S.\ Sunder, whose
answer proved rather more difficult than the author originally imagined.
The author thanks them for the question, for their patience and for their encouragement.

\section{Amalgamated free product of a II$_1$--factor and another von Neumann algebra}

The (reduced) amalgamated free product construction for C$^*$--algebras and conditional expectations 
appeared in the earliest work~\cite{V85}
of Voiculescu on freeness;
taking the closure in the strong--operator topology, one gets the amalgamated free product of von Neumann algebras.
See~\cite{V95} for more, or the book~\cite{VDN}.
Some considerations involving normality of conditional expectations, embeddings and
completely positive maps of amalgamted free products of von Neumann algebras are treated in~\cite{BlD01}.

Amalgamated free products of finite von Neumann algebras with respect to trace--preserving
conditional expectations were used in work of Popa~\cite{Po93} as part of a construction of subfactors.
We will use the expression
\begin{equation}\label{eq:M=ADB}
\Mcal=A*_DB
\end{equation}
to mean
$A$ and $B$ are von Neumann algebras with faithful, normal, tracial states $\tau_A$ and $\tau_B$,
respectively and
$D$ is von Neumann algebra that is embedded in both $A$ and $B$ as
unital von Neumann subalgebras so that $\tau_A$ and $\tau_B$ restrict to the same trace on $D$;
we let $E^A_D:A\to D$ and $E^B_D:B\to D$ be the trace--preserving conditional expectations onto $D$
and we take the amalgmated free product of von Neumann algebras
\[
(\Mcal,E_D)=(A,E^A_D)*_D(B,E^B_D).
\]
Then $E_D$ is a conditional expectation from $\Mcal$ onto $D$ and $\tau:=\tau_A\circ E_D$
is a faithful, normal, tracial state on $\Mcal$.
In fact, we always take $D$ to be finite dimensional.

If $D=\Cpx$ in~\eqref{eq:M=ADB}, then we say $\Mcal$ is a free product over the scalars,
and we may write simply $\Mcal=A*B$.
It follows from results of~\cite{D94} that $\Mcal=A*B$ is a II$_1$--factor if $B$ is a II$_1$--factor and $A$ is any finite
von Neumann algebra.

\begin{thm}\label{thm:M=ADB}
Let $\Mcal=A*_DB$ where $B$ is a II$_1$--factor and $D$ is finite dimensional.
Then $\Mcal$ is a II$_1$--factor.
Moreover, $\Mcal$ can be described recursively
in terms of the operations of compression/dilation and taking free products
over the scalars, by using the following facts:
\begin{enumerate}[(a)]
\item
If $q$ is an abelian projection in $D$ whose central carrier in $D$ is $1$, then
\[
q\Mcal q=(qAq)*_{qDq}(qBq),
\]
so we may without loss of generality assume $D$ is commutative,
and we do assume this in items (\ref{item2})--(\ref{itemlast}) below.
Furthermore, we assume $\dim D>1$.

\item \label{item2}
Suppose $p$ is a minimal projection of $D$.
Let
\[
\Nc_1=W^*((1-p)A(1-p)\cup(1-p)B(1-p)).
\]
Then
\[
\Nc_1\cong (1-p)A(1-p)*_{(1-p)D}\,(1-p)B(1-p)
\]
and $\Nc_1$ is a II$_1$--factor.

\item
Suppose there is another minimal projection $p_0$ of $D$ with $\tau(p_0)\ge\tau(p)$.
Let $r=C_A(1-p)p$ be the part of the central carrier in $A$ of the projection $1-p$ that lies under $p$.
Let
\[
\At=(1-p)A(1-p)+\Cpx r +(p-r)A(p-r)
\]
and let $\Nc_2=W^*(\At\cup B)$.
Let $y\in B$ be a partial isometry such that $y^*y=p$ and $f:=yy^*\le p_0$.
Then $\Nc_2$ is a II$_1$--factor and
\begin{equation}\label{eq:N2fp}
f\Nc_2f\cong f\Nc_1f*p\At p,
\end{equation}
where $p\At p$ is endowed with the tracial state obtained by restricting $\tau_A$.

\item \label{itemlast}
If $r\ne0$, then
\begin{equation}\label{eq:rMr}
r\Mcal r\cong r\Nc_2 r*L(\Ints).
\end{equation}
\end{enumerate}
\end{thm}

Before embarking on the proof, we describe some technical notation.

\begin{notation}
For subsets $X_1$ and $X_2$ of an algebra, we use the notation $\Lambdao(X_1,X_2)$
for the set of all words of the form $c=c_1c_2\ldots c_n$ for $n\in\Nats=\{1,2,\ldots\}$,
and $c_j\in X_{i(j)}$ for some $i(1),\ldots,i(n)\in\{1,2\}$ with $i(j)\ne i(j+1)$.
Sometimes such a product will be thought of as an element of the algebra and sometimes
as a formal word, with letters $c_1,\ldots,c_n$.

For any subalgebra $Q\subseteq\Mcal$ we write $Q\ominus D$ for $Q\cap\ker E_D$ and $Q\oup$ for $Q\cap\ker\tau$.
\end{notation}

\begin{proof}[Proof of Theorem~\ref{thm:M=ADB}]
Part~(a) is certainly well known ({\em cf}\/ Lemma 5.2 of~\cite{D95}), but we'll run through an argument.
It is clear that $qAq$ and $qBq$ are free with amalgamation over $qDq$ with respect to the restiction of $E_D$ to $q\Mcal q$.
The set 
\begin{equation}\label{eq:LoAB}
D+\Lambdao(A\ominus D,B\ominus D)
\end{equation}
is a strong--operator--topology (s.o.t.)\ dense
$*$--subalgebra of $\Mcal$.
There is a finite list of partial isometries $v_0,v_1,\ldots,v_\ell\in D$
such that $v_j^*v_j\le q$, $v_0=q$ and $\sum_{j=1}^\ell v_jv_v^*=1$.
Using these, we easily write every element of the set~\eqref{eq:LoAB} as a linear combination of elements
from
\[
\bigcup_{i,j=0}^\ell v_i(qDq+\Lambdao(qAq\ominus D,qBq\ominus D))v_j^*.
\]
From this, we see that $W^*(qAq\cup qBq)=q\Mcal q$ holds, and~(a) is proved.

The proof of the rest of this theorem (for commutative $D$) proceeds by induction on the dimension of $D$.
Then the assertion of part~(b) is just the induction hypothesis.
However, to initiate the induction argument, we need to treat the case $D=\Cpx$ and show that a free product
$A*B$ is a II$_1$--factor when $B$ is a II$_1$--factor.
As already mentioned, this follows from results of~\cite{D94}.

For part~(c), we have $B=W^*((1-p)B(1-p)\cup\{y\})$, so
$\Nc_2=W^*(\Nc_1\cup\{y\}\cup p\At p)$
and from this we obtain $f\Nc_2f=W^*(f\Nc_1f\cup y(p\At p)y^*)$.
For~\eqref{eq:N2fp},
it will suffice to show that $f\Nc_1f$ and $y\At y^*$ are free with respect to
$\tau(p)^{-1}\tau\restrict_{f\Mcal f}$, which is equivalent to showing that we have
\begin{equation}\label{eq:yfN1}
\Lambdao(y^*(f\Nc_1f)\oup y,(p\At p)\oup)\subseteq\ker\tau.
\end{equation}
However, since $f$ is a subprojection
of a minimal projection in $D$, we have
\[
E_D((f\Nc_1f)\oup)=0
\]
and
(using Kaplansky's density theorem) every element of $(f\Nc_1f)\oup$ is the s.o.t.\ limit of a bounded
sequence in
\[
\lspan\Lambdao((1-p)A(1-p)\ominus D,(1-p)B(1-p)\ominus D).
\]
But $E_D(y^*(f\Nc_1f)\oup y)=0$, so $y^*E_B((f\Nc_1 f)\oup)y=E_B(y^*(f\Nc_1 f)\oup y)\subseteq B\ominus D$,
and $(1-p)By\subseteq B\ominus D$.
Thus, we see that every element of $y^*(f\Nc_1f)\oup y$ is the s.o.t.\ limit of a bounded
sequence in $\lspan\Theta$, where $\Theta$ is the set of words in
\[
\Lambdao((1-p)A(1-p)\ominus D,B\ominus D)
\]
whose first and last letters are from $B\ominus D$.
So, in order to show~\eqref{eq:yfN1}, it will suffice to show
$\Lambdao(\Theta,(p\At p)\oup)\subseteq\ker\tau$.
But since $(p\At p)\oup\subseteq A\ominus D$, we have
\[
\Lambdao(\Theta,(p\At p)\oup)\subseteq\Lambdao(\Theta,A\ominus D)\subseteq\Lambdao(A\ominus D,B\ominus D)\subseteq\ker\tau,
\]
and~\eqref{eq:N2fp} is proved.
Now $f\Nc_2f$ is a factor by results of~\cite{D94};
since $f$ is full in $\Nc_2$, it follows that $\Nc_2$ is a factor, and~(c) is proved.

To prove~(d), first note that there are partial isometries $v_i\in A$ ($i\in I$), for $I$ a finite or countable index set,
such that $r=\sum_{i\in I}v_i^*v_i$ and $v_iv_i^*\le 1-p$, and then $A=W^*(\At\cup\{v_i\mid i\in I\})$.
Since $\tau(p)\le1/2$ and $\Nc_2$ is a factor, we may choose $x_i\in\Nc_2$ such that $x_ix_i^*=v_iv_i^*$ and the projections
$(x_i^*x_i)_{i\in I}$ are pairwise orthogonal.
Since $r\in\Nc_2$ and $\sum_{i\in I}x_i^*x_i$ is a projection in $\Nc_2$ having the same trace as $r$,
we can choose $y\in\Nc_2$ such that $y^*y=r$ and $yy^*=\sum_{i\in I}x_i^*x_i$.
Let $w=\sum_{i\in I}y^*x_i^*v_i$.
Then $w^*w=ww^*=r$ and, moreover,
$\Mcal=W^*(\Nc_2\cup\{w\})$.
We will show that, in $r\Mcal r$ with respect to the trace $\tau(r)^{-1}\tau\restrict_{r\Mcal r}$,
the element
$w$ is a Haar unitary and is $*$--free from $r\Nc_2 r$.
This will imply~\eqref{eq:rMr}.
So it will suffice to show
\[
\Lambdao((r\Nc_2r)\oup,\{w^n\mid n\in\Nats\}\cup\{(w^*)^n\mid n\in\Nats\})\subseteq\ker\tau.
\]
For this, it will suffice to show $\Theta\subseteq\ker\tau$,
where $\Theta$ is the set of all words in
\[
\Lambdao(r\Nc_2r,\{w,w^*\})
\]
such that
every letter from $r\Nc_2r$ that has $w$ on the left and $w^*$ on the right, or $w^*$ on the left and $w$ on the right,
belongs to $(r\Nc_2r)\oup$.

Take a word from the set $\Theta$ described above.
Now write out $w=\sum_iy^*x_i^*v_i$ and $w^*=\sum_iv_i^*x_iy$ for every $w$ and $w^*$ except those that appear as part of a sequence
$w^*zw$, where $z\in(r\Nc_2r)\oup$.
So doing, we realize the word as equal to a sum of (or, in the case of infinite $I$,
the s.o.t.\ limit of a bounded sequence consisting of sums of) words from
\begin{equation}\label{eq:bigLam}
\Lambdao\big((r\Nc_2r)\oup\cup\bigcup_{i\in I}x_iy\Nc_2r\cup\bigcup_{i\in I}r\Nc_2y^*x_i^*\,,
\,\{w,w^*\}\cup\{v_i\mid i\in I\}\cup\{v_i^*\mid i\in I\}\big)
\end{equation}
that satisfy
\begin{enumerate}[(i)]
\item each letter $w$ or $w^*$ can appear only as part of a sequence $w^*zw$ for $z\in(r\Nc_2r)\oup$,
\item if a letter comes from $(r\Nc_2r)\oup$ and is not part of a sequence $w^*zw$ as in~(i),
then it (a) has $w$ or a $v_i$ on the left or is the first letter
and (b) has $w^*$ or a $v_j^*$ on the right or is the last letter
\item if a letter comes from $r\Nc_2y^*x_j^*$, then it (a) has a $w$ or a $v_i$ on the left or is the first letter
and (b) has $v_j$ on the right
\item if a letter comes from $x_iy\Nc_2r$, then it (a) has $v_i^*$ on the left and (b) has a $w^*$ or a $v_j^*$ on the right or is the last letter
\item no letter has a $v_i^*$ on the left and a $v_j$ on the right.
\end{enumerate}
It will, therefore, suffice to show that each word from~\eqref{eq:bigLam} that satisfies (i)-(v) above
evaluates to zero under $\tau$.

We now describe certain approximations of elements satisfying the various conditions above.

We treat $(r\Nc_2r)\oup$ from~(ii).
Since $r$ is a minimal projection of $\At$, we have
\[
E_A((r\Nc_2r)\oup)=E_{\At}((r\Nc_2r)\oup)=\{0\}.
\]
Therefore, every element of $(r\Nc_2r)\oup$ is the s.o.t.\ limit of a bounded sequence in
$\lspan(\Lambdao(\At\ominus D,B\ominus D)\backslash \At\ominus D)$.
Moreover, since $r\At=\Cpx r$, we may left-- and right--multiply by $r$ to see that every element of $(r\Nc_2r)\oup$
is the s.o.t.\ limit of a bounded sequence in $\lspan\Bc$, where $\Bc$ is the set of words in $\Lambdao(\At\ominus D,B\ominus D)$
whose first and last letters belong to $B\ominus D$.

We treat $r\Nc_2y^*x_j^*$ from~(iii).
Since $r\At(1-p)=0$ and $x_jx_j^*\le 1-p$, we find $E_A(r\Nc_2y^*x_j^*)=E_{\At}(r\Nc_2y^*x_j^*)=0$ and
then, by left--multiplying with $r$, we find that every element of $r\Nc_2y^*x_j^*$ is the s.o.t.\ limit
of a bounded sequence in $\lspan\Bc_1$, where $\Bc_1$ is the set of words in
$\Lambdao(\At\ominus D,B\ominus D)$
whose first letter belongs to $B\ominus D$.

We similarly treat $x_iy\Nc_2r$ from~(iv).
Every element of $x_iy\Nc_2 r$ is the s.o.t.\ limit
of a bounded sequence in $\lspan\Bc_2$, where $\Bc_2$ is the set of words in
$\Lambdao(\At\ominus D,B\ominus D)$
whose last letter belongs to $B\ominus D$.

We treat $w^*zw$ from~(i).
Consider
\[
w^*zw=\sum_{i,j\in I}v_i^*x_iyzy^*x_j^*v_j\,,
\]
for $z\in(r\Nc_2r)\oup$.
We have $\tau(w^*zw)=0$;
since $r$ is a subprojection of a minimal projection in $D$ and $w=rwr$, we have $E_D(w^*zw)=0$.
So we have
\begin{equation}\label{eq:Ewzw}
E_A(w^*zw)=\sum_{i,j}v_i^*E_{\At}(x_iyzy^*x_j^*)v_j\in A\ominus D.
\end{equation}
For each $i$ and $j$, the element $x_iyzy^*x_j^*-E_{\At}(x_iyzy^*x_j^*)\in\Nc_2$, is the s.o.t.\ limit
of a bounded sequence in $\lspan(\Lambdao(\At\ominus D,B\ominus D)\backslash \At\ominus D)$.
Because $r$ is a subprojection of the minimal projection $p$ of $D$,
we have $\At v_j=\At(1-p)v_jp=(1-p)\At v_jp\subseteq A\ominus D$ and
$v_i^*\At\subseteq A\ominus D$.
Therefore, $v_i^*(x_iyzy^*x_j^*-E_{\At}(x_iyzy^*x_j^*))v_j$ is the s.o.t.\ limit of a bounded sequence
in $\lspan\Ac$, where $\Ac$ is the set of words in $\Lambdao(A\ominus D,B\ominus D)$ whose first and last letters
are from $A\ominus D$.
We would like to conclude that
\[
w^*zw-E_A(w^*zw)=\sum_{i,j\in I}v_i^*(x_iyzy^*x_j^*-E_{\At}(x_iyzy^*x_j^*))v_j
\]
is the s.o.t.\ limit of a bounded sequence in $\lspan\Ac$.
However, as $I$ may be infinite, there is at first sight difficulty with the boundedness.
This can be overcome as follows:
every element of $\Mcal$ is the s.o.t.\ limit of a bounded sequence in $\lspan(D\cup\Lambdao(A\ominus D,B\ominus D))$.
Since each $v_i^*(x_iyzy^*x_j^*-E_{\At}(x_iyzy^*x_j^*))v_j$ is the s.o.t.\ limit of a bounded sequence in $\lspan\Ac$, it follows
that it is orthogonal to $D$ and to all elements of $\Lambdao(A\ominus D,B\ominus D)$ that are not in $\Ac$.
From this, it follows that $w^*zw-E_A(w^*zw)$ is the s.o.t.\ limit of a bounded sequence on $\lspan\Ac$.
However, since $A\ominus D\subseteq\Ac$, from~\eqref{eq:Ewzw} we have that $w^*zw$ is itself the s.o.t.\ limit
of a bounded sequence in $\lspan\Ac$.

\smallskip
Let $c$ be a word from~\eqref{eq:bigLam} that satisfies (i)--(v).
We will show $\tau(c)=0$.
In light of the treatments of the several terms above, it will suffice to show $\tau(\ct)=0$,
whenever $\ct$ is obtained from $c$ by replacing
\begin{enumerate}[{r}-i.]
\item any sequence of letters $w^*zw$ for $z\in(r\Nc_2r)\oup$ with an arbitrary element of $\Ac$,
\item any letter from $(r\Nc_2r)\oup$ that in $c$ has a $w$ or a $v_i$ on the left or is the first letter and has a $w^*$ or a $v_j^*$ on the right or is the last letter,
with an arbitrary element of $\Bc$
\item any letter from $r\Nc_2yx_j$ that in $c$ has a $w$ or a $v_i$ on the left or is the first letter and has $v_j^*$ on the right,
with an arbitrary element of $\Bc_1$,
\item any letter from $x_i^*y^*\Nc_2r$ that in $c$ has $v_i$ on the left and $w^*$ or a $v_j^*$ on the right or is the last letter
with an arbitrary element of $\Bc_2$.
\end{enumerate}
Then $\ct$ can be seen as an element of
\[
\Lambdao(\Ac\cup\{v_i\mid i\in \}\cup\{v_j^*\mid j\in I\}\,,\,\Bc_1\cup\Bc_2)
\]
(we note that $\Bc=\Bc_1\cap\Bc_2$)
such that every letter from $\Bc_1\backslash\Bc$ is followed by a $v_i$ and every letter from $\Bc_2\backslash\Bc$
is preceded by a $v_j^*$.
Now using again $\At v_i\subseteq A\ominus D$ and $v_j^*\At\subseteq A\ominus D$, we see that the word $\ct$ can
be seen as an element of $\Lambdao(A\ominus D,B\ominus D)$, and, therefore, $\tau(\ct)=0$.
\end{proof}

\section{Old results and ``generating sets of free dimension''}
\label{sec:fdim}

In this section, we recall some old results and the ``free dimension'' assignment.

The interpolated free group factors $L(\Fb_t)$, ($1<t\le\infty$) form a family of 
II$_1$--factors related to each other by compression/dilation, via the formula
\[
L(\Fb_t)_s=L(\Fb_{1+s^{-2}(t-1)}),\qquad (1<t\le\infty,\,0<s<\infty),
\]
and such that if $t=n$ is an integer, then $L(\Fb_t)$ is the usual group von Neumann algebra of the
nonabelian free group on $n$ generators.
(See~\cite{Ra94} and~\cite{D94:i}).
It is known that either $L(\Fb_r)\cong L(\Fb_t)$ for all $r,t\in(1,\infty]$ or $L(\Fb_r)\not\cong L(\Fb_t)$
for all $1<r<t\le\infty$.
(This was proved by R\u adulescu in~\cite{Ra94}, and in slightly weaker form in~\cite{D94:i}.)

\begin{notation}\label{not:S}
Let $\SEu$ be the class of pairs $(A,\tau)$ of finite von Neumann algebras $A$
and normal, faithful, tracial states $\tau$, such that $A$ is a direct sum of finitely or countably infinitely many von Neumann algebras, each of which is one of the following:
\begin{enumerate}[(i)]
\item finite dimensional
\item a diffuse hyperfinite von Neumann algebra
\item an interpolated free group factor.
\end{enumerate}
\end{notation}

\begin{thm}[\cite{D93}]
\label{thm:D93}
Let $(A,\tau_A)$ and $(B,\tau_B)$ belong to the class $\SEu$ with $\dim(A)\ge2$ and $\dim(B)\ge3$
and let
\[
(\Mcal,\tau)=(A,\tau_A)*(B,\tau_B)
\]
be the free product.
Then $\Mcal=L(\Fb_t)\oplus D$,
where $D$ is either $0$ or is finite dimensional.
\end{thm}
Moreover, in~\cite{D93}, an algorithm is given to find $D$ and the restriction of $\tau$ to $D$.
Also, an algorithm is given to find the parameter $t$, in terms of the initial data.
This latter algorithm is most easily described by a quanitity
that was, perhaps misleadingly, called ``free dimension,'' denoted $\fdim$,
and was additive for free products:
$\fdim(\Mcal,\tau)=\fdim(A,\tau_A)+\fdim(B,\fdim_B)$.
The problem is that we assigned $\fdim(L(\Fb_t))=t$, which may be nonsense.
But because the only purpose of free dimension was to decide what the parameter $t$ is in an interpolated free group factor  $L(\Fb_t)$, 
if it's nonsense, then it's mathematically harmless nonsense.
Nonetheless, it understandably caused some confusion, and was in that sense not harmless.

The solution proposed and implemented in~\cite{D02} for this (at the suggestion of a referee) was to replace
the symbols $\fdim(\Mcal,\tau)=s$ with the words ``$\Mcal$ (with its trace $\tau$) has a generating set of
free dimension $s$.''
This convention, though more cumbersome, is mathematically defensible, because the only generating sets to which we assign a ``free dimension'' are
those obtained by certain operations (free products, rescalings, direct sums,
certain sorts of inductive limits, certain limits in strong--operator--topology,
{\em etc.}), for which we have proved results that yield isomorphisms.
See the discussion in section~1 of~\cite{D02} for more on this;
in particular, see the rules (i)--(v) for possession of ``generating sets of free dimension'' found on p.\ 147 of~\cite{D02}.
We would like to point out the following rule, which follows directly from rule~(iv), cited above:
\begin{itemize}
\item[(iv')]
if $p$ is a central projection of $\Mcal$, with $\tau(p)=\alpha$, and if $p\Mcal$ and $(1-p)\Mcal$
(each endowed with the appropriate renormalization of the restriction of $\tau$)
have generating sets of free dimension $d_1$ and $d_2$, respectively,
then $(\Mcal,\tau)$ has a generating set of free dimension $1+\alpha^2(d_1-1)+(1-\alpha)^2(d_2-1)$.
\end{itemize}
It will be helpful to have the following rule for rescaling, which will be familiar
to those accustomed to interpolated free group factors;
the proof serves as a quick summary of the rules.
\begin{lemma}\label{lem:fdimcutdown}
Suppose $(A,\tau_A)\in\SEu$
has a generating set of free dimension $t$
and let $p$ be a projection in $A$ which central support is $1$, and with $\tau_A(p)=\beta$.
Then $(pAp,\beta^{-1}\tau_A\restrict_B)$ has a generating set of free dimension $1+\beta^{-2}(t-1)$.
\end{lemma}
\begin{proof}
If $A$ is an interpolated free group factor, then this follows from the usual formula for rescaling them.
If $A$ is a matrix algebra $M_n(\Cpx)$ (and thus has free dimension $1-n^{-2}$)
or a diffuse hyperfinite von Neumann algebra (and thus has free dimension $1$),
then this is an easy calculation.
If $A$ is a finite or countably infinite direct sum
\[
A=\bigoplus_{i\in I}\smd{A_i}{\alpha_i},
\]
where each $A_i$ is a matrix algebra or a diffuse hyperfinite von Neumann algebra or an interpolated free group factor,
then each $A_i$ has generating set of some free dimension $t_i$, and $t=1+\sum_i\alpha_{i\in I}^2(t_i-1)$.
We write $p=\bigoplus_{i\in I}p_i$ and let $\beta_i=\tau_A(p_i)$, so that $\sum_i\beta_i=\beta$ and
\[
pAp=\bigoplus_{i\in I}\smd{p_iA_ip_i}{\beta_i/\beta}.
\]
As observed above, each $p_iA_ip_i$ has a generating set of free dimension $1+\big(\frac{\alpha_i}{\beta_i}\big)^2(t_i-1)$.
Thus, $pAp$ has a generating set of free dimension
\[
1+\sum_{i\in I}\big(\frac{\beta_i}\beta)^2\big(\big(\frac{\alpha_i}{\beta_i}\big)^2(t_i-1)\big)
=1+\sum_{i\in I}\big(\frac{\alpha_i}\beta)^2(t_i-1)=1+\beta^{-2}(t-1).
\]
\end{proof}

\begin{remark}
Kenley Jung showed~\cite{J03} that for a hyperfinite algebra $A\in\SEu$ with trace $\tau_A$, all finite generating sets
have the same free entropy dimension of Voiculescu (see~\cite{V02}),
and the value of this free entropy dimension coincides with the ``free dimension''
assigned in the above--described scheme.
Thus, for such algebras, we use the notation $\fdim(A,\tau_A)$, or simply $\fdim(A)$ when the trace
$\tau_A$ is clearly indicated by the context, for this value of the free entropy dimension.
\end{remark}

\section{Amalgamated free products of certain direct sums}

The result of the algorithm described in Theorem~\ref{thm:M=ADB} can be described precisely when the component algebras $A$ and $B$
are built up out of hyperfinite von Neumann algebras, by using the results recalled in Section~\ref{sec:fdim} and
the rules for generating sets of free dimension.
The class $\SEu$ is described above in Notation~\ref{not:S}.

\begin{prop}\label{prop:M=ADB}
Let $\Mcal=A*_DB$ where $D$ is finite dimensional,
where $B$ is either an interpolated free group factor or the hyperfinite II$_1$--factor
and where $(A,\tau_A)\in\SEu$, $A\ne D$.
Suppose $(A,\tau_A)$ has a generating set of free dimension $x$.
Then $\Mcal$ is an interpolated free group factor $L(\Fb_t)$, where letting $\tau_D$ denote the restriction of $\tau_A$ to $D$, we have
\[
t=\begin{cases}
1+x-\fdim(D,\tau_D),&\text{if }B\text{ is hyperfinite} \\
s+x-\fdim(D,\tau_D),&\text{if }B=L(\Fb_s).
\end{cases}
\]
\end{prop}
\begin{proof}
Part~(a) of Theorem~\ref{thm:M=ADB} together with Lemma~\ref{lem:fdimcutdown} shows that we may without loss
of generality assume $D$ is commutative.
If $D=\Cpx$, then this was proved in~\cite{D93} (see Theorem~\ref{thm:D93} above).
For higher dimensional commutative $D$, we proceed by induction on the dimension of $D$,
using the other parts of Theorem~\ref{thm:M=ADB}.

For the remainder of this proof we will use the convention that $L(\Fb_1)$ means the hyperfinite II$_1$--factor,
and we write $B=L(\Fb_s)$ for $1\le s\le\infty$.
We write
\begin{equation}\label{eq:A}
A=\smdp{A_0}{\alpha_0}{q_0}\oplus\bigoplus_{i\in I}\smdp{A_i}{\alpha_i}{q_i}
\end{equation}
where $A_0$ is diffuse hyperfinite (or zero), $I$ is finite or countably infinite, $0\not\in I$ and for each $i\in I$,
$A_i$ is either a matrix algebra $M_{n_i}(\Cpx)$,
in which case we set $x_i=\fdim(A_i)=1-n_i^{-2}$,
or an interpolated free group factor $L(\Fb_{x_i})$.
Thus, $A_i$ has a generating set of free dimension $x_i$.
The notation of~\eqref{eq:A} indicates that $q_i$ is the central projection of $A$ that is the identity element of the summand $A_i$
and $\tau(q_i)=\alpha_i$.
Then $A$ has a generating set of free dimension $x$, where
\[
x=1+\sum_{i\in I}\alpha_i^2(x_i-1).
\]

Let $p$ be a minimal projection of $D$ of smallest trace value, and write $p=p_0+\sum_{i\in I}p_i$ where $p_i=pq_i$.
Let $\beta_i=\tau(p_i)$ and $\beta=\tau(p)$.
We will use part~(b) of Theorem~\ref{thm:M=ADB} and the induction hypothesis to find $\Nc_1$, as defined in that theorem.
We easily compute
\[
\fdim\big((1-p)D\big)=1+(1-\beta)^{-2}(\fdim(D)-1+\beta^2).
\]
and $(1-p)B(1-p)=L(\Fb_{1+(1-\beta)^{-2}(s-1)})$.
Let $I'=\{i\in I\mid p_i\ne q_i\}$.
Then
\[
(1-p)A(1-p)=(q_0-p_0)A_0(q_0-p_0)\oplus\bigoplus_{i\in I'}(q_i-p_i)A_i(q_i-p_i)
\]
and $(1-p)A(1-p)$ has a generating set of free dimenstion
\[
1+(1-\beta)^{-2}\sum_{i\in I'}\alpha_i^2(x_i-1).
\]
Therefore, $\Nc_1=L(\Fb_{t_1})$, where
\[
t_1=1+(1-\beta)^{-2}\sum_{i\in I'}\alpha_i^2(x_i-1)+(1-\beta)^{-2}(s-1)-(1-\beta)^{-2}(\fdim(D)-1+\beta^2).
\]

We now use part~(c) of Theorem~\ref{thm:M=ADB} to find $\Nc_2$, as defined there.
Let $\ell=C_{A_0}(q_0-p_0)p_0$ be the part of the central carrier of $q_0-p_0$ that lies under $p_0$.
Then $r=\ell+\sum_{i\in I'}p_i$.
Let $\gamma=\tau(r)$.
We find
\[
p\At p=\smdp{\Cpx}{\gamma/\beta}{r}\oplus(p_0-\ell)A_0(p_0-\ell)
\oplus\bigoplus_{i\in I\backslash I'}\smd{p_iA_ip_i}{\beta_i/\beta}.
\]
and we calculate that $p\At p$ has a generating set of free dimension
\[
1-(\gamma/\beta)^2+\beta^{-2}\sum_{i\in I\backslash I'}\alpha_i^2(x_i-1).
\]
So by the initial step of Theorem~2.1 (i.e.\ when $D=\Cpx$) and the rescaling formula for interpolated free
group factors, we find
$f\Nc_2f=L(\Fb_{t_2})$, where
\[
t_2=1-(\gamma/\beta)^2+\beta^{-2}(s-1)+\beta^{-2}(x-1)-\beta^{-2}(\fdim(D)-1).
\]

Now further rescaling and part~(d) of Theorem~\ref{thm:M=ADB} finish the proof.
\end{proof}

\begin{lemma}\label{lem:minproj}
Let $B$, $C_1$ and $C_2$ be finite von Neumann algebras, let $k\in\Nats$ and let
\begin{align*}
A&=C_1\oplus M_k(\Cpx)\otimes C_2 \\
A_0&=C_1\oplus M_k(\Cpx),
\end{align*}
where we regard $A_0$ as a subalgebra of $A$.
Suppose a von Neumann algebra $D$ is embedded as a unital subalgebra of both $A_0$ and $B$ and suppose
faithful traces on $A$ and $B$ are specified that restrict to the same trace on $D$.
Let $\Mcal_0=A_0*_DB$ and $\Mcal=A*_DB$ be the amalgmated free products.
Let $p$ be a minimal (nonzero) projection in $0\oplus M_k(\Cpx)\subset A_0$

Then $p\Mcal p\cong p\Mcal_0p*C_2$.
More specifically, $p\Mcal p$ is generated by $p\Mcal_0 p$ and $p\otimes C_2$, and these two algebras are
free with respect to $\tau(p)^{-1}\tau\restrict_{p\Mcal p}$.
\end{lemma}
\begin{proof}
Clearly, $p$ is equivalent in $A_0$ to a subprojection of some minimal projection in $D$,
so we may without loss of generality assume
that $p\le q$ for a minimal projection $q$ of $D$.
We must show
\[
\Lambdao\big((p\Mcal_0p)\oup,p\otimes C_2\oup\big)\subseteq\ker\tau.
\]
Since $p$ is mininal in $A_0$,
every element of $(p\Mcal_0p)\oup$ is the s.o.t.\ limit of a bounded sequence is
$\lspan\Theta$,
where $\Theta$ is the set of words in $\Lambdao(A_0\ominus D,B\ominus D)$ whose first and last letters are
from $B\ominus D$.
So it will suffice to show
$\Lambdao\big(\Theta,p\otimes C_2\oup\big)\subseteq\ker\tau$.
But since $p\otimes C_2\subseteq A\ominus D$, we have
\[
\Lambdao\big(\Theta,p\otimes C_2\oup\big)\subseteq\Lambdao(A\ominus D,B\ominus D)\subseteq\ker\tau.
\]
\end{proof}

\begin{notation}
Let $\REu$ be the class of pairs $(A,\tau)$ of finite von Neumann algebras $A$
and normal, faithful, tracial states $\tau$, such that $A$ is a direct sum of finitely many von Neumann algebras, each of which is
one of the following:
\begin{enumerate}[(i)]
\item finite dimensional
\item $L^\infty([0,1])\otimes M_n(\Cpx)$ for some $n\in\Nats$
\item the hyperfinite II$_1$--factor
\item an interpolated free group factor.
\end{enumerate}
Let $\REu_0$ be the class of pairs $(A,\tau)$ in $\REu$ but where all interpolated free group factors $L(\Fb_{t_i})$
appearing as summands in $A$ have finite parameter $t_i$.
\end{notation}

\begin{thm}\label{thm:R0}
The class $\REu$ is closed under taking free products with amalgamation over finite dimensional algebras.
Moreover, if $(A,\tau_A)$ and $(B,\tau_B)$ belong to $\REu$ and have generating sets of free dimension
$d_A$ and $d_B$, respectively, and if
$\Mcal=A*_DB$ is an amalgamted free product over a
finite dimensional subalgebra $D$, then $\Mcal$ has a generating set of free dimension
\begin{equation}\label{eq:fdim}
d_A+d_B-\fdim(D).
\end{equation}
Finally, also $\REu_0$ is closed under taking free products with amalgamation over finite dimensional algebras.
\end{thm}
\begin{proof}
Let $\Mcal=A*_DB$ for $A$ and $B$ in the class $\REu$, and as usual denote the trace resulting from the free product construction
by $\tau$.
If $A$ and $B$ are finite dimensional, then the facts that $(\Mcal,\tau)$
lies in $\REu_0$ and has a generating set of the correct free dimension were proved in~\cite{D95}.

Suppose $A$ and $B$ in $\REu_0$ are both direct sums of finite dimensional and finite type~I algebras, i.e., of those
in classes~(i) and~(ii), above.
Then $A$ can be written as an inductive limit $A=\overline{\bigcup_{k\ge1}A_k}$ of an increasing sequence
$A_1\subseteq A_2\subseteq\cdots$
of finite dimensional subalgebras, with $D\subseteq A_1$, and where each inclusion $A_k\hookrightarrow A_{k+1}$
is of the form
\begin{multline}\label{eq:fdincl}
A_k=M_{\ell(1)}(\Cpx)\,\oplus\, M_{\ell(2)}(\Cpx)\,\oplus\,\cdots\,\oplus\, M_{\ell(p)}(\Cpx) \\
\hookrightarrow
\big(M_{\ell(1)}(\Cpx)\otimes C\big)\,\oplus\, M_{\ell(2)}(\Cpx)\,\oplus\,\cdots\,\oplus\, M_{\ell(p)}(\Cpx)=A_{k+1},
\end{multline}
where $C=\Cpx\oplus\Cpx$ is a two--dimensional algebra with minimal projections having traces $1/2$.
Also $B$ can be similarly written.
The free product $A*_DB$ is, thus, the inductive limit of the amalgamated free products $\Mcal_k=A_k*_DB_k$
of finite dimensional algebras.
As we saw above, each algebra $\Mcal_k$ is in the class $\REu_0$, and has a generating set of free dimension
$\fdim(A_k)+\fdim(B_k)-\fdim(D)$.
Moreover, for $k$ large enough, the finite dimensional part (if any), the minimal central projections and the 
type I part (of any) of $\Mcal_k$ are fixed,
and only parameters $t$ in the interpolated free group factor summands $L(\Fb_t)$ of $\Mcal_k$ can vary with $k$.
We will see that these changes in $t$ are from so--called standard embeddings (see~\cite{D93}).
If $p$ is a minimal projection in the direct summand $M_{\ell(1)}(\Cpx)$ of $A_k$ in~\eqref{eq:fdincl} above,
then, by Lemma~\ref{lem:minproj}, the subalgebra $p(A_k*_DB_k)p\subseteq p(A_{k+1}*_DB_k)p$
is freely complemented by a copy of $\Cpx\oplus\Cpx$.
By Proposition~4.4 of~\cite{D93}, this inclusion is a standard embedding.
Therefore, the dilation of this inclusion to the central carrier $q=C_{\Mcal_k}(p)$, namely
$q(A_k*_DB_k)q\hookrightarrow q(A_{k+1}*_DB_k)q$ is an standard embedding (see Proposition~4.2 of~\cite{D93})
and also the inclusion $q\Mcal_kq\hookrightarrow q\Mcal_{k+1}q$ is an standard embedding.
Since the direct limit of interpolated free group factors $L(\Fb_{t_k})$ under standard embeddings is the interpolated free group factor,
$L(\Fb_t)$ where $t$ is the limit of the increasing sequence $t_k$,
(see Proposition~4.3 of~\cite{D93}), we have that $\Mcal$ is in $\REu$.
Since $\lim_{k\to\infty}\fdim(A_k)=\fdim(A)$, and similarly for $B$, $\Mcal$ has a generating set of free dimension $\fdim(A)+\fdim(B)-\fdim(D)$.

In the general case, we proceed by induction on the total number of summands in $A$ and $B$ that are hyperfinite
II$_1$--factors or interpolated free group factors.
The case of no such summands was treated above.
Suppose there is at least one, and without loss of generality suppose $A$ has one.
Let $p$ be a central projection of $A$ so that $pA$ is an interpolated free group factor or hyperfinite II$_1$--factor.
By Lemma~4.3 of~\cite{D95}, 
\begin{equation}\label{eq:pp}
p\Mcal p\cong pA*_{pD}(p\Mcalt p),
\end{equation}
where $\Mcalt=(pD\oplus(1-p)A)*_DB$
By the induction hypothesis, $\Mcalt$ (with its trace) belongs to $\REu$, and Proposition~\ref{prop:M=ADB} implies that $p\Mcal p$ is
an interpolated free group factor.
Of course, the central support $C_\Mcal(p)$ of $p$ in $\Mcal$ equals the central support of $p$ in $\Mcalt$, and
$\Mcal\cong\Mcalt(1-C_\Mcal(p))\oplus C_\Mcal(p)\Mcal$, where the last summand is an interpolated free group factor.
This implies that $\Mcal$ is in the class $\REu$.

It remains to see that the ``free dimension'' calculations work out, and that $\Mcal$ belongs to $\REu_0$ if $A$ and $B$ do.
Let $\alpha=\tau(p)$ and suppose $A$ and $B$ have generating sets of free dimensions $d_A$ and $d_B$, respectively.
By assumption, $pA=L(\Fb_s)$ for some $s\in[1,\infty]$, where (as in the proof of Prop.~\ref{prop:M=ADB}),
we use the convention that $L(\Fb_1)$ denotes the hyperfinite II$_1$--factor.
Let $\At=pD\oplus(1-p)A$;
then by rule~(iv'), $\At$ has a generating set of free dimension $d_A+\alpha^2(\fdim(pD)-s)$.
Therefore, by the induction hypothesis, $\Mcalt$ has a generating set of free dimension
\begin{equation}\label{eq:f1}
d_{\Mcalt}:=d_A+\alpha^2(\fdim(pD)-s)+d_B-\fdim(D).
\end{equation}
Let $q=C_{\Mcal}(p)$ and let $\beta=\tau(q)$.
Then, by rule~(iv'), 
\begin{equation}\label{eq:f2}
d_{\Mcalt}=1+\beta^2(d_{q\Mcalt}-1)+(1-\beta)^2(d_{(1-q)\Mcalt}-1),
\end{equation}
where $q\Mcalt$ and $(1-q)\Mcalt$ have generating sets of free $d_{q\Mcalt}$ and $d_{(1-q)\Mcalt}$, respectively.
By Lemma~\ref{lem:fdimcutdown}, $p\Mcalt p$ has a generating set of free dimension
$d_{p\Mcalt p}:=1+(\beta/\alpha)^2(d_{q\Mcalt}-1)$.
By~\eqref{eq:pp} and Proposition~\ref{prop:M=ADB}, $p\Mcal p$ is the interpolated free group factor with parameter
\begin{equation}\label{eq:f3}
s+d_{p\Mcalt p}-\fdim(pD)=s+1+(\beta/\alpha)^2(d_{q\Mcalt}-1)-\fdim(pD)
\end{equation}
and, consequently, $q\Mcal$ is the interpolated free group factor with parameter
\begin{equation}\label{eq:f4}
d_{q\Mcal}:=(\alpha/\beta)^2(s-\fdim(pD))+d_{q\Mcalt}.
\end{equation}
Using $(1-q)\Mcal=(1-q)\Mcalt$ and rule~(iv'), we find that $\Mcal$ has a generating set of free dimension
\[
d_{\Mcal}:=1+\beta^2(d_{q\Mcal}-1)+(1-\beta)^2(d_{(1-q)\Mcalt}-1)=d_A+d_B-\fdim(D),
\]
as required,
where the last equality was obtained using~\eqref{eq:f1}--\eqref{eq:f4}.
Finally, if we started with $A$ and $B$ in $\REu_0$, then all of the values obtained in~\eqref{eq:f3} and~\eqref{eq:f4}
for parameters of interpolated free group factors are finite, so the free product algebra $\Mcal$ is also in $\REu_0$.
\end{proof}

\begin{remark}
The formula~\eqref{eq:fdim}
is consistent with the main result of~\cite{BDJ08}.
Indeed, the later paper shows that if $\Mcal=A*_DB$ is a tracial amalgamated free product
of von Neumann algebras with $D$ hyperfinite and if $X$, $Y$ and $Z$ are finite generating sets for $A$, $B$ and $D$, respectively,
then
\[
\delta_0(X\cup Y\cup Z)=\delta_0(X\cup Z)+\delta_0(Y\cup Z)-\delta_0(Z).
\]
All von Neumann algebras in the family $\REu_0$ posses finite generating sets $X$ whose free entropy
dimensions are the expected ``free dimension'' used in Theorem~\ref{thm:R0}.
However, the main content of Theorem~\ref{thm:R0}, namely,
that the free product is in the class $\REu$ (or $\REu_0$), does not follow from the free dimension calculation alone.

Indeed, an example of N.\ Brown~\cite{B05} shows that there is a II$_1$--factor $\Mcal$ having finite generating set $X$
with free entropy dimension $1<\delta_0(X)<\infty$, and even having finite free entropy, but so that $\Mcal$ is not an interpolated
free group factor.
Of course, it is also entirely possible that there is an interpolated free group factor $L(\Fb_t)$ having a finite generating set $X$
with $\delta_0(X)\ne t$.
In principle, this could occur even if the free group factors are not isomorphic to each other.
\end{remark}

\end{document}